\documentclass[12pt, oneside,reqno]{amsart}
\usepackage{hyperref}
\usepackage{geometry}
\usepackage{amsmath}
\usepackage{cite}
\newcommand{\rn}{\mathbb R^n}
\newcommand{\sn}{S^{n-1}}
\newcommand{\kno}{\mathcal K^n_o}
\newcommand{\kne}{\mathcal K^n_e}
\newcommand{\conv}{\operatorname{conv}}

\newcommand{\hm}{\mathcal H^{n-1}}
\newcommand{\rbo}{\mathbb R}
\newcommand{\balpha}{\pmb{\alpha}}

\newcommand\wtilde[1]{\overset{\lower.4ex\hbox{$\scriptstyle \sim$}}{#1}}

\numberwithin{equation}{section}

\newtheorem{thm}{Theorem}[section]

\newtheorem{lem}[thm]{Lemma}

\newtheorem{defi}[thm]{Definition}

\begin{document}
\title[Gaussian measure and $L_p$-Minkowski problem  ]{The $L_p$-Gaussian Minkowski problem }
\author[J. Liu]{Jiaqian Liu}
\address{School of Mathematics, Hunan University, Changsha, 410082, Hunan Province, China}
\email{liujiaqian@hnu.edu.cn}

\keywords {Minkowski problem, $L_p$-Gaussian Minkowski problem, Monge-Amp\`ere equation, degree theory}
\subjclass[2010] {52A40, 52A38, 35J96.}

\thanks{The author was partially supported by Tian Yuan Special Foundation(12026412) and Hunan Science and Technology Planning Project(2019RS3016).}

\maketitle
\begin{abstract}
In this paper,  we extend the article that Minkowski problem in Gaussian probability space of Huang et al.  to $L_p$-Gaussian Minkowski problem, and obtain the existence and uniqueness of $o$-symmetry weak solution  in case of $p\geq1$.
\end{abstract}
\section{Introduction}
The Brunn-Minkowski theory mainly studies the geometric functional of convex bodies and their differential. The differentiable functional will generate some  intrinsic geometric measures. In Euclidean space $\rn$, the most significant functional of convex body is volume denoted by $V$,
 and Aleksandrov \cite{A38} gave the  variational formula
\begin{equation}
\lim_{t\rightarrow 0}\frac{V(K+tL)-V(K)}{t}=\int_{\sn}h_L(v)dS_K(v),
\end{equation}
where $K + tL = \{x + ty :\text{$x\in K$ and $y \in L$}\}$, $h_L:\sn\rightarrow \rbo$ is the support function of convex body $L$. It  naturally generates the surface area measure $S_K$ of convex body $K$. The classical Minkowski problem characterizes the surface area measure, it reads : {\it given a finite Borel measure $\mu$ on $\sn$, what are the necessary and sufficient conditions on $\mu$ does there exist  a convex body $K$ in $\rn$ such that surface area measure $S_K(\cdot\,)=\mu(\cdot\,)$? If $K$ exists, is it unique?}   In the smooth case, it is equivalent to solve a certain fully nonlinear elliptic equation. The classical Minkowski problem was solved by Minkowski  in \cite{HM}  for polytope, and later solved by Aleksandrov \cite{A38},  Fenchel and Jessen \cite{FJ38}  in general.

In 1993, Lutwak developed the Brunn-Minkowski theory to $L_p$-Brunn-Minkowski theory. The $L_p$-Minkowski problem, characterizing $L_p$-surface area measure $S_p(K)={h_K}^{1-p}S_K$. Clearly, when $p=1$, the  $L_p$-surface area measure is just the classical surface measure.  Lutwak \cite{Lut1} testified that there exists a unique $o$-symmetry convex body such that  $S_p(K)=\mu$ for $p>1$. The $L_p$-Minkowski problem has been studied extensively in  \cite{BHZ16,B-LYZ,CLZ17, HLX, HLYZ,JLW,SA, Zhu,Zhu14,Zhu17} and so on. The uniqueness of the Minkowski problem is related to Brunn-Minkowski inequality. For $p>1$, the $L_p$-Brunn-Minkowski inequality can be obtained easily. In case of $0<p<1$, $L_p$-Brunn-Minkowski inequality and its condition of the equality holds have not been improved much until  \cite{CHL}. While the logarithmic Minkowski problem $(p = 0)$ and the centro-affine Minkowski problem $(p =- n)$ are two
special  cases; see, e.g., \cite {BLYZ13}, and \cite {CW}. The regularity of the $L_p$-Minkowski problem, for example \cite{CW, HQL, LV}.

In \cite{Lut}, the dual Brunn-Minkowski theory was developed in the 1970s. The most significant  dual curvature measure and its associated Minkowski problem  in the dual Brunn-Minkowski theory were studied  in \cite{HYYZ}.  The work of \cite{HYYZ} play a key role in promoting the development of partial differential equation and dual
Brunn-Minkowski theory as \cite{CHZ,HZ,BLYZ19,LSW,CL,WFZ}. Recently,  a number of meaningful and interesting work  emerged on the  Orlicz-Minkowski problem, such as  \cite{GHXY1,GHXY2,HCLYZ,JL,LL,XY,ZXY}.

In \cite{HXZ}, the authors posed the Minkowski problem in Gaussian  probability space. Gaussian volume functional denoted  by $\gamma_n$ plays a central role, which is given by
\[
\gamma_n(K)=\frac{1}{(\sqrt{2\pi})^n}\int_K e^{\frac{-|x|^2}{2}}dx.
\]
Using the variational formula from Huang et al. \cite{HYYZ}, the variational formula of $\gamma_n$ is obtained  as follows:
\begin{equation}
\lim_{t\rightarrow 0}\frac{\gamma_n(K+tL)-\gamma_n(K)}{t}=\int_{\sn}h_L(v)dS_{\gamma_n,K}(v),
\end{equation}
 where Gaussian surface area measure $S_{\gamma_n,K}$ is given by
 $$
S_{\gamma_n,K}(\eta)=\frac{1}{(\sqrt{2\pi})^n}\int_{\nu_K^{-1}(\eta)}e^{\frac{-|x|^2}{2}}d\hm(x).
$$

About the mainly result  of Gaussian Minkowski problem showed  in \cite{HXZ} Theorem 1.4, stated as: Let $\mu$ be a finite even Borel measure on $\sn$ that is not concentrated on any closed hemisphere and $|\mu|<\frac{1}{\sqrt{2\pi}}$. Then there exists a unique $o$-symmetry convex body with $\gamma_n(K)>\frac{1}{2}$ such that $S_{\gamma_n,K}=\mu.$

As remarked in \cite{HXZ},  Gaussian Minkowski problem is different from the Minkowski problem in Lebesgue measure space in these aspects. At first,  the Gaussian surface area of all convex set in $\rn$ is not more than $4{n^{\frac{1}{4}}}$ which is shown detailed in Ball \cite{B} and \cite{F}. That is to say, the allowable $\mu$ in Gaussian Minkowski problem can not have an arbitrarily big total mass. Secondly, Gaussian probability measure has neither translation invariance nor homogeneity, which makes our research more difficult. As a final remark, by the definition of Gaussian surface area, we can know that not only small convex bodies can have smaller Gaussian surface areas, but also large convex bodies. For instance, let $B_r$ be a centered ball with radius $r,$ the Gaussian surface area's density  of $B_r$ is $f_r=\frac{1}{(\sqrt{2\pi})^n}e^{\frac{-r^2}{2}}r^{n-1}$. One can observe that $\frac{1}{(\sqrt{2\pi})^n}e^{\frac{-r^2}{2}}r^{n-1}\rightarrow 0$ as $r\rightarrow 0 $ or $r\rightarrow \infty.$

Naturally, we develop the $L_p$-Brunn-Minkowski theory in Gaussian probability space. The purpose of this paper is to continue the study in \cite{HXZ} and to construct $L_p$-Gaussian surface area measures.
A  $L_p$-variational formula  similar to Lemma 3.2 in Lutwak \cite{Lut1} is gained as below:
$$
\lim _{t \rightarrow 0} \frac{\gamma_n\left([h_t] \right)-\gamma_n\left(K\right)}{t}= \frac{1}{p}\int_{\sn}f(v)^p dS_{p,\gamma_n,K}(v),
$$
where $h_t(v)=(h_K(v)^p+tf(v)^p)^{\frac{1}{p}}$ is $L_p$-combination for $p\neq0$, and  it leads to define the $L_p$-Gaussian surface area measures as follows
$$
S_{p,\gamma_n,K}(\eta)=\frac{1}{(\sqrt{2\pi})^n}\int_{\nu_K^{-1}(\eta)}e^{\frac{-|x|^2}{2}}(x\cdot\nu_K(x))^{1-p}d\hm(x).
$$

The $L_p$-Gaussian Minkowski problem  reads : {\it given a finite Borel measure $\mu$ on $\sn$, what are the necessary and sufficient conditions on $\mu$ does there exist  a convex body $K$ in $\rn$ such that $L_p$-Gaussian surface area measure $S_{p,\gamma_n,K}(\cdot\,)=\mu(\cdot\,)$? If $K$ exist, is it unique?}

The first assertion  involves the uniqueness result, it is known to us that the Gaussian surface area's density $f_r$ of $B_r$ tends to $0$ as $r\rightarrow 0$ or $r\rightarrow \infty,$ this means  that Gaussian surface area is may not unique. Our study in the restricted context that the Gaussian volume grater than or equal to $\frac{1}{2}$, then obtain the following uniqueness result. The same fact was firstly found in \cite{HXZ} in the situation of $p=1$.

\begin{thm}[\bf Uniqueness] For $p\geq1$, let $K,L \in \kno$  with the same $L_p$-Gaussian surface area, i.e.,
\begin{equation*}
S_{p,\gamma_n,K}=S_{p,\gamma_n,L}.
\end{equation*}
If $\gamma_n(K), \gamma_n(L)\geq\frac{1}{2},$  then  $K=L.$
\end{thm}
Aleksandrov's variational method implies the existence of solution with lagrange multipliers.
\begin{thm}{\it For $ p>0,$ let $\mu$ be a nonzero finite  Borel measure on $\sn$ and be  not concentrated in any closed hemisphere. Then there exist a  $K\in\kno$  and a positive constant $\lambda$ such that
\begin{equation*}
\mu=\frac{\lambda}{p} S_{p,\gamma_n,K}.
\end{equation*}
}
\end{thm}
Due to the lack of homogeneity of measure, just like all previous Orlicz-Minkowski problems, in the solving process, the lagrange multipliers produced by Aleksandrov's variational method  can not eliminated. Recently, the work of Huang et al. in \cite{HXZ},  the first one overcame this obstruct. Inspired their work,   we obtain the existence of smooth solution.
\begin{thm}[\bf Existence of smooth solution]
For $p\geq1,$  let even function $f\in C^{2,\alpha}(\sn)$ is positive  and $|f|_{L_1}<\sqrt{\frac{2}{\pi}}r^{-p}ae^{\frac{-a^2}{2}}$, the $r$ and $a$ are chosen such that $\gamma_n(rB)=\gamma_n(P)=\frac{1}{2}$, symmetry strip $P=\left\{x\in \rn:|x_1|\leq a\right\}$. Then there exists a unique $ C^{4,\alpha}$ convex body $ K\in \kne $ with $\gamma_n(K)>\frac{1}{2}$, the support function $h_K$ solves the equation
\begin{equation}
\frac{1}{(\sqrt{2\pi})^n}{h_K}^{1-p}e^{\frac{-(|\nabla h_K|^2+h_K^2)}{2}}\det(\nabla^2h_K+h_KI)=f.
\end{equation}
\end{thm}
Finally, the existence of weak solution of the $L_p$-Gaussian Minkowski problem follows by an approximation method with the help of a uniform estimate.
\begin{thm}[\bf Existence of weak solution ]{\it For $p\geq1,$ let $\mu$ be a nonzero finite even  Borel measure on $\sn$ and be not concentrated in any closed hemisphere with $|\mu|<\sqrt{\frac{2}{\pi}}r^{-p}ae^{\frac{-a^2}{2}},$ the $r$ and $a$ are chosen such that $\gamma_n(rB)=\gamma_n(P)=\frac{1}{2}$, symmetry strip $P=\left\{x\in \rn:|x_1|\leq a\right\}$.} Then there exists a unique $K\in\kne$ with $\gamma_n(K)>\frac{1}{2}$ such that
\begin{equation*}
S_{p,\gamma_n,K}=\mu.
\end{equation*}
\end{thm}

We  overcome the most important obstacle of  the $L_p$ Gaussian  Minkowski problem  about whether the solution degenerates or not, and obtain that it holds for asymmetric convex bodies see Theorem \ref{Eoo}. But  the limitation of this item is the lack of $L_p$  Gaussian isoperimetric type inequality for general convex body such that it can only achieve the $o$-symmetry convex body, see Theorem \ref{13} and Theorem \ref{exi}.

 The paper  is organized  into five sections: we list some conceptual knowledge of convex bodies in  Preliminaries part,  and gain the key  $L_p$-variational formula  in Section~\ref{SecP}. In Section~\ref{SecMP},  Aleksandrov's variational method is applied to obtain the non-symmetry solution with lagrange multipliers for $p>0$. The uniqueness of solution is proved in Section~\ref{SecU1}. Section~\ref{SecE} provides the existence result of $o$-symmetry smooth solution for $p\geq1$ by means of the degree theory of partial differential equation, then  the $o$-symmetry weak solution is obtained with the help of approximation argument.

%\section{Preliminaries}

\section{$L_p$-Gaussian surface area and the $L_p$-Gaussian Minkowski problem }

\label{SecP}
In this section, we make a brief review of some relevant notions about convex bodies.    More detailed information of convex body theory can be found in the  Schneider \cite{Sch}.

Let $\rn$ represent the $n$-dimensional Euclidean space, $\sn$ be the unit sphere in $\rn,$ $C^+(\sn)$ be positive functions defined on $\sn$ and $|\mu|$ denote total measure of  the Borel measure $\mu$ on $\sn.$

Convex body in $\rn$ is  compact convex set with nonempty interior denoted by $\mathcal K^n$. Let $\partial K$ denote the boundary of $K$, $\kno$ be convex bodies that contain the origin in interior, and $\kne$ denote the $o$-symmetry convex bodies.

The support function  $h_K$ of a compact convex subset $K$ in $\rn$ is defined as
\[
h_K(y)=\max_{x\in K} y\cdot x
\]
for $y \in\rn$, where $\cdot$ represents the inner product.

 For $K\in\kno$, the radial function $\rho_K$ is defined by
\[
\rho_K(x) = \max\{\lambda : \lambda x \in K\}
\]
for $x\in \rn\setminus\{0\}$. A simple observation allows us to obtain
\[
\partial K =\{\rho_K(u)u : u\in \sn\}.
\]

In fact, $\mathcal K^n$ is the space of compact convex set in $\rn$ and it can be equipped with {\it Hausdorff metric} which associated norm is
 \[
\delta(K,L) =\max_{v\in\sn}|h_K(v)-h_L(v)|.
\]
  We can easily show that $\delta(K,L)$ is the distance between $K,L\subset\mathcal K^n$.
  Therefore, the convergence of compact convex set  sequence can be characterized by the convergence of support function. That is to say,  $K_i\in\mathcal K^n\rightarrow K\in\mathcal K^n$ if
\begin{equation*}
\max_{v\in\sn}|h_{K_{i}}(v)-h_{K}(v)| \rightarrow 0,
\end{equation*}
as $i\rightarrow\infty.$ Or if  $K_i\in \kno$,
\begin{equation*}
\max_{u\in\sn}|\rho_{K_{i}}(u)-\rho_{K}(u)| \rightarrow 0,
\end{equation*}
as $i\rightarrow\infty,$ then $K_i\rightarrow K $.

Let $h\in C^+(\sn)$, the {\it Wulff shape } $[h]$ generated by $h$   is a convex body defined by

\[
[h]= \bigcap\limits_{v\in\sn}\left\{x\in\rn: x\cdot v\leq h(v)\right\}.
\]
Obviously, if $K\in\kno$,
\[
[ h_K]= K.
\]

Let $\rho\in C^+(\sn)$, the {\it convex hull } $\langle\rho\rangle$  generated by $\rho$ is a convex body defined by
\[
\langle\rho \rangle  = \conv\{\rho(u)u : u\in\sn\}.
\]
Clearly, if $K\in\kno$,
\begin{equation*}
\langle \rho_K \rangle = K.
\end{equation*}

For $K, L\subset \mathcal K^n$ and real $a,b>0$, the {\it Minkowski combination}, $aK+bL\subset \mathcal K^n$ defined as follow
\[
aK + bL = \{ax + by :\text{$x\in K$ and $y \in L$}\},
\]
and support function is
\begin{equation*}
h_{aK+bL} = ah_K + bh_L.
\end{equation*}

Firey was the first to define the  {\it $L_p$-combination} of  $K,L\in\kno$, for $a,b>0,$ $p\neq0,$
\begin{equation*}
a\cdot K+_pb\cdot L=\bigcap\limits_{v\in\sn}\left\{x\in\rn: x\cdot v\leq \left(a{h_K(v)}^p+b{h_L(v)}^p\right)^{\frac{1}{p}}\right\},
\end{equation*}
and
\begin{equation*}
a\cdot K+_0b\cdot L=\bigcap\limits_{v\in\sn}\left\{x\in\rn: x\cdot v\leq {h_K(v)}^a{h_L(v)}^b\right\}.
\end{equation*}
Here $a\cdot K=a^{\frac{1}{p}}K$.

In case of $p\geq 1,$ $\left(a{h_K(v)}^p+b{h_L(v)}^p\right)^{\frac{1}{p}}$ is just the support function, but for $0<p<1$ is not necessarily.

For  $K \in \kno$, its
{\it polar body} $K^* \in \kno$  in $\rn$ is defined by
\[
K^* = \{ x\in \rn : x\cdot y \le 1, \  \text{for all }\  y\in K \}.
\]

It is easy to show that, if $K\in\kno$,
\begin{equation}\label{polar-identity}
 \rho_K = 1/h_{K^*}\quad\text{and}\quad h_K= 1/\rho_{K^*}
\end{equation}
for $x\in\rn\setminus\{0\}$. In addition, we can get  $(K^*)^*=K$ if $K\in\kno$.

For $K, K_i\in\kno$,  by \eqref{polar-identity} and the previous characterization of convex body convergence, then  $K_i\rightarrow K $ if and only if $K_i^*\rightarrow K^*. $

For convex body  $K\in\rn$, then its {\it supporting hyperplane  with outward unit normal vector $v\in \sn$} represented by
\begin{equation*}
H_K(v) = \{x \in \rn : x\cdot v = h_K(v)\}.
\end{equation*}

  The boundary point  of $K$  only has one  supporting hyperplane called regularity point, otherwise, it is a singular point.  The set of singular points denoted as $\sigma_K$, it is well known that $\sigma_K$ has spherical Lebesgue measure $0$ (see Schneider \cite{Sch}, p.84).

  For  $x\in {\partial K\setminus\sigma_K } $, its {\it Gauss map}  $\nu_K$  is represented by
 \begin{equation*}
\nu_K(x) =\{v \in \sn: x\cdot v = h_K(v) \}.
 \end{equation*}

 Correspondingly, for Borel set $\eta\subset\sn$, its {\it reverse Gauss map} denoted by $\nu_K^{-1}$,
 \begin{equation*}
 \nu_K^{-1}(\eta)=\{x \in \partial K: \nu_K(x)\in\eta \}.
 \end{equation*}

 For Borel set $\eta\subset\sn$, its surface area measure is defined as
\begin{equation*}
S_K(\eta) =\mathcal H^{n-1}(\nu_K^{-1}(\eta)),
\end{equation*}
where $\mathcal H^{n-1}$ is $(n-1)$-dimensional Hausdorff measure.

  Then $g\in C(\sn),$
 \begin{equation}\label{normal}
\int_{\partial K\setminus\sigma_K} g(\nu_K(x))\, d\hm(x) = \int_{\sn} g(v)\,dS_K(v).
\end{equation}

For Borel set $\omega\subset\sn$,  $\balpha_K(\omega)$  denotes its {\it radial Gauss image }  and  is defined as
\begin{equation*}
\balpha_K(\omega) = \{v\in \sn : \rho_K(u)(u\cdot v) = h_K(v) \},
\end{equation*}
 for  $u\in\omega$.  When Borel set $\omega$ have only one element $u$, we will  abbreviate $\balpha_{K}(\{u\})$ as $\balpha_K(u)$. The subset of $\sn$ which  make $\balpha_K(u)$ contain more than one element denoted by $\omega_K$ for each $u\in\omega$. The set $\omega_K$ has spherical Lebesgue measure $0$ (see Schneider \cite{Sch}, Theorem 2.2.5).

The {\it radial Guass map} of $K$ is a map denoted by $\alpha_K(u)$, the only difference between $\alpha_K(u)$ and $\balpha_K(u)$ is that the former is defined on $\sn\setminus \omega_K$ not on $\sn$ which lead to $\balpha_K(u)$ may have many elements but $\alpha_K(u)$ has only one. In other words, if $\balpha_K(u)=(\{v\})$, then $\balpha_K(u)=\alpha_K(u)$.

For $\hm$-integrable function $g: \partial K\rightarrow \rbo$,
\begin{equation}\label{radial}
\int_{\partial K}g(x)d\hm(x)=\int_{\sn}g(\rho_K(u)u)F(u)du,
\end{equation}
where $F$ is defined $\hm$-a.e. on $\sn$ by
\begin{equation*}
F(u)=\frac{(\rho_K(u))^n}{h_K(\alpha_K(u))}.
\end{equation*}

For Borel set $\eta\subset\sn$, its {\it reverse radial Gauss image }  $\balpha^*_K(\eta)$  is shown as
\begin{equation*}
\balpha^*_K(\eta) = \{u\in \sn :\rho_K(u) (u\cdot v) = h_K(v) \}
\end{equation*}
for $v\in\eta$.

When Borel set $\eta$ have only one element $v$, we will  abbreviate $\balpha^*_K(\{v\})$ as $\balpha^*_K(v)$. The subset of $\sn$ which  make $\balpha^*_K(v)$ contain more than one element denoted by $\eta_K$ for each $v\in\eta$. The spherical Lebesgue measure of set $\eta_K$ is  $0$ (see Schneider \cite{Sch}, Theorem 2.2.11).

The {\it reverse radial Gauss map} of $K$ is a map denoted by $\alpha^*_K(v)$,  if $\balpha^*_K(v)=(\{u\})$, then $\balpha^*_K(v)=\alpha^*_K(v)$.

%\section{$L_p$-Gaussian surface area and the $L_p$-Gaussian Minkowski problem }\label{SecV}
For this part, we give the definition of the $L_p$-Gaussian surface area measure as follows.
\begin{defi}
Let $K\in\kno$, for Borel set $\eta\in\sn$, $p\in\mathbb R$, define $L_p$-Gaussian  surface area as below
\begin{equation*}
S_{p,\gamma_n,K}(\eta)=\frac{1}{(\sqrt{2\pi})^n}\int_{\nu_K^{-1}(\eta)} e^{-\frac{|x|^2}{2}}(x\cdot\nu_K(x))^{1-p}d\hm(x),
\end{equation*}
which is a Borel measure on $\sn$.
\end{defi}

The above definition is deduced from the $L_p$-variational formula, and the proof of variational formula requires the following lemma. Our proof mainly relies on ideas developed in \cite{HYYZ}.
\begin{lem}\label{variation}
For $p\neq0,$ let $K\in\kno$, $f: \sn \rightarrow \mathbb R$ be a continuous function. For enough small $\delta>0$, and each $t\in(-\delta,\delta)$,  define the continuous function $h_t: \sn\rightarrow(0,\infty)$ as
\begin{equation*}
h_t(v)=(h_K(v)^p+tf(v)^p)^{\frac{1}{p}}\quad v\in\sn.
\end{equation*}
Then,
\begin{equation*}
\lim _{t \rightarrow 0} \frac{\rho_{[h_t] }(u)-\rho_K(u)}{t}=\frac{f(\alpha_K(u))^p}{ph_K(\alpha_K(u))^p}\rho_K(u)
\end{equation*}
holds for almost all $u\in\sn$. In addition, there exists $M>0$ such that
\begin{equation*}
|\rho_{[h_t]}(u)-\rho_K(u)|\leq M|t|,
\end{equation*}
for all $u\in\sn$ and $t\in(-\delta,\delta)$.
\end{lem}
\begin{proof}
Since
\begin{equation*}
h_t(v)=(h_K(v)^p+tf(v)^p)^{\frac{1}{p}}
\end{equation*}
is equivalent to
\begin{equation*}
\log h_t(v)=\log h_K(v)+\frac{tf(v)^p}{p{h_K(v)}^p},
\end{equation*}
 the outcome is an easy application of the Lemmas 2.8 and 4.1 in \cite{HYYZ}.
\end{proof}
The following $L_p$-variational formula, as we shall see in the next section, is the key to solve the $L_p$-Gaussian Minkowski problem.
\begin{thm}[\bf $L_p$-variational formula in Gaussian space ]\label{Variation}
For $p\neq0,$ let $K\in\kno$, $f: \sn \rightarrow \mathbb R$ be continuous function. For sufficiently small $\delta>0$, and each $t\in(-\delta,\delta)$,  define the continuous function $h_t: \sn\rightarrow(0,\infty)$ by
\begin{equation*}
h_t(v)=(h_K(v)^p+tf(v)^p)^{\frac{1}{p}}\quad v\in\sn.
\end{equation*}
 Then,
\begin{equation}\label{4VFE}
\lim _{t \rightarrow 0} \frac{\gamma_n\left([h_t] \right)-\gamma_n\left(K\right)}{t}= \frac{1}{p}\int_{\sn}f(v)^p dS_{p,\gamma_n,K}(v).
\end{equation}
\end{thm}
\begin{proof}
Applying polar coordinates, shows that
\begin{equation*}
\gamma_n([h_t])=\frac{1}{(\sqrt{2\pi})^n}\int_{\sn}\int_{0}^{\rho_{[h_t]}(u)}e^{\frac{-r^2}{2}}r^{n-1}drdu.
\end{equation*}
Since $K\in\kno$ and $h_t\rightarrow h_K$ uniformly as $t\rightarrow 0$, by {\it Aleksandrov's Convergence Lemma}, we can get $[h_t]\rightarrow [h_K]=K$, then radial function of $[h_t]$ converges to the radial function of $K$, i.e., $\rho_{[h_t]}\to \rho_K$,  and there exist $m_0, m_1>0$ such that  $\rho_{[h_t]}, \rho_K\in[m_0,m_1]$.

For simplify, denote $F(s)=\int_{0}^{s}e^{\frac{-r^2}{2}}r^{n-1}dr$. By mean value theorem, there exists $\theta\in[m_0,m_1]$ such that
\begin{equation*}
|F(\rho_{[h_t]}(u))-F(\rho_K(u))|= |F^{'}(\theta)||\rho_{[h_t]}(u)-\rho_K(u)|\leq M|F^{'}(\theta)||t|,
\end{equation*}
 the $M$ is arisen as in Lemma \ref{variation}. With the simple calculation,
 $$
 |F^{'}(\theta)|=|e^{\frac{-\theta^2}{2}}\theta^{n-1}|\leq M _1
 $$
 for some  constant $M_1>0$ as $\theta\in[m_0,m_1]$.

  Therefore,
\begin{equation*}
|F(\rho_{[h_t]}(u))-F(\rho_K(u))|\leq MM_1|t|.
\end{equation*}
 Together with Lemma \ref{variation}, employing the dominated convergence theorem, it follows  that
 \begin{equation*}
 \begin{split}
 \lim _{t \rightarrow 0} \frac{\gamma_n\left([h_t] \right)-\gamma_n\left(K\right)}{t}&=\frac{1}{(\sqrt{2\pi})^n}\int_{\sn}\frac{f(\alpha_K(u))^p}{ph_K(\alpha_K(u))^p}\rho_K(u)^ne^{\frac{-\rho_K(u)^2}{2}}du\\
 &=\frac{1}{p(\sqrt{2\pi})^n}\int_{\partial K}f(\nu_K(x))^pe^{\frac{-|x|^2}{2}}(x\cdot\nu_K(x))^{1-p}d\hm(x)\\
 &=\frac{1}{p}\int_{\sn}f(v)^pdS_{p,\gamma_n,K}(v).
 \end{split}
 \end{equation*}
\end{proof}

By the definition of $S_{p,\gamma_n,K}$, together with the conclusions that the Gauss surface area measure is weakly convergent and is absolutely continuous which were got in \cite{HXZ}, it sufficient to reveal that $S_{p,\gamma_n,K}$ is weakly convergent with respect Hausdorff metric, and is absolutely continuous with respect to surface area measure.
\section{The variational method to generate solution}
\label{SecMP}
This section aims to transform the $L_p$-Gaussian Minkowski problem into an optimization  problem by employing variational method, and prove  the optimizer is just the solution to the Minkowski problem of $S_{p,\gamma_n,K}$.
\subsection{An associated optimization problem}
For any nonzero finite  Borel measure $\mu$ on  $\sn$, $p>0,$  define $\phi: \kno\rightarrow \mathbb{R}$
by
\begin{equation*}
\phi(Q)=\int_{\sn}h_Q(v)^pd\mu(v),
\end{equation*}
for each $Q\in\kno $.
We take into account the following minimum problem
\begin{equation}
\min\left\{\phi(Q) : \text{$ Q \in \kno$ }, \gamma_n(Q) =\frac{1}{2}\right\}.
\end{equation}
\begin{lem}\label{If}
{\it If $K\in\kno$  and satisfies
 \begin{equation*}
 \phi(K)=\min\left\{\phi(Q) : \text{$ Q \in \kno$ }, \gamma_n(Q) =\frac{1}{2} \right\},
 \end{equation*}
equivalent to
 \begin{equation*}
\varphi(h_K)=\min\left\{\varphi(z) : \text{$ z \in C^{+}(\sn)$ }, \gamma_n([z]) =\frac{1}{2}\right \},
\end{equation*}
where $\varphi:C^{+}(\sn)\rightarrow \mathbb{R}$  is defined by
\begin{equation*}
\varphi(z)=\int_{\sn}z(v)^pd\mu(v).
\end{equation*}
}
\end{lem}
\begin{proof}
For the Wulff shape
\[
[z]= \bigcap\limits_{v\in\sn}\left\{x\in\rn: x\cdot v\leq z(v)\right\},
\]
it is easy to  get $h_{[z]}\leq z$ and $[h_{[z]}]=[z]$, thus we have
\begin{equation*}
\varphi(z)\geq\varphi(h_{[z]}).
\end{equation*}
%this conclusion implies us that
%\begin{equation*}
%\inf\{\Psi_f(z) : \text{$ z \in C_e^{+}(\sn)$ and $\int_{[z]} f(x)dx=1$} \}=\inf\{\Psi_f(h_Q) : \text{$ Q \in \mathcal{K}_{e}^{n}$ and $\int_Q f(x)dx=1$} \}.
%\end{equation*}

Clearly, $K\in\kno$  and satisfies
 \begin{equation*}
\phi(K)=\min\left\{\phi(Q) : \text{$ Q \in \kno$ }, \gamma_n(Q) =\frac{1}{2} \right \}
 \end{equation*}
 if and only if
 \begin{equation*}
\varphi(h_K)=\min\left\{\varphi(z) : \text{$ z \in C^{+}(\sn)$ }, \gamma_n([z]) =\frac{1}{2}  \right\}.
\end{equation*}
\end{proof}
\begin{lem}\label{var}
{\it For $p>0,$ let $\mu$ be a nonzero finite  Borel measure on $\sn$.  If $K\in\kno$  and satisfies
 \begin{equation}\label{bi}
 \phi(K)=\min\left\{\phi(Q) : \text{$ Q \in \kno$ }, \gamma_n(Q) =\frac{1}{2} \right \},
 \end{equation}
 then there exists a constant $\lambda>0$   such that $\mu=\frac{\lambda S_{p,\gamma_n,K}}{p}$}.
\end{lem}
\begin{proof}
By \eqref{bi} and Lemma \ref{If}, we  get
\begin{equation}\label{max}
\varphi(h_K)=\min\left\{\varphi(z) : \text{$ z \in C^{+}(\sn)$ }, \gamma_n([z]) =\frac{1}{2} \right \}.
\end{equation}

For any $f\in C(\sn)$ and $t\in(-\delta,\delta)$ where $\delta>0$ is sufficiently small, let
\begin{equation*}
 h_t(v)=(h_K(v)^p+tf(v)^p)^{\frac{1}{p}}.
 \end{equation*}

 \eqref{max} implies that $h_K$ being a minimizer,   there exists a constant $\lambda>0$ such that
 \begin{equation*}
\left. \frac{d}{d t} \right|_{t=0}\varphi(h_t)=\lambda\left. \frac{d}{d t} \right|_{t=0}\gamma_n([h_t]).
  \end{equation*}
Moreover,  from the $L_p$-variational formula \eqref{4VFE}, we have
   $$
   \int_{\sn}f(v)^pd\mu(v)=\frac{\lambda}{p}\int_{\sn}f(v)^pdS_{p,\gamma_n,K}(v).
   $$

  By the arbitrariness of $f$, then  $\mu=\frac{\lambda}{p} S_{p,\gamma_n,K}$.
 \end{proof}

\subsection{Existence of an optimizer }
\begin{thm}\label{Eoo}
 For $ p>0,$  let $\mu$ be a nonzero finite  Borel measure on $\sn$ and be not concentrated on any closed hemisphere. Then there exists a $K\in\kno$ such that
\[
\phi(K)=\min\left\{\phi(Q) : \text{$ Q \in \kno$ }, \gamma_n(Q) =\frac{1}{2}\right \}.
\]
\end{thm}
\begin{proof}
Suppose ${Q_l}\subset\kno$ is a minimal sequence, i.e.,
\begin{equation}\label{ms}
\lim_{l\rightarrow\infty}\phi(Q_l)=\min\left\{\phi(Q) : \text{$ Q \in \kno$ }, \gamma_n(Q) =\frac{1}{2}\right \}.
\end{equation}

We  now claim that $Q_l$ is uniformly bounded. If not, then there exists $u_l\in\sn$
such that $\rho_{Q_l}(u_l)\rightarrow\infty$ as $l\rightarrow\infty$. By the definition of support function,  $\rho_{Q_l}(u_l)(u_l\cdot v)_+\leq h_{Q_l}(v)$, then
 \begin{equation*}
 \begin{split}
\phi(Q_l)&=\int_{\sn}h_{Q_l}(v)^pd\mu(v)\\
&\geq\int_{\sn}\rho_{Q_l}(u_l)^p{(u_l\cdot v)_+}^pd\mu(v)\\
&= \rho_{Q_l}(u_l)^p\int_{\sn}{(u_l\cdot v)_+}^pd\mu(v),
\end{split}
 \end{equation*}
 where $(t)_+=\max\{{t,0}$\} for any $t\in \mathbb{R}$. Since  $\mu$ is not concentrated on any closed hemisphere, there exists a positive constant $c_0$ such that
  $$
 \int_{\sn}{(u_l\cdot v)_+}^pd\mu(v)>c_0.
 $$

 Therefore,
 \begin{equation*}
\phi(Q_l)>\rho_{Q_l}(u_l)^pc_0\rightarrow \infty
 \end{equation*}
 as $l\rightarrow\infty.$  But this is contradicted to \eqref{ms}.
Then we conclude $Q_l$ is uniformly bounded. By Blaschke selection theorem, $Q_l$ has convergent subsequence,  still denoted by $Q_l$,  converges to  a compact convex set $K$ of $\rn$.
By the continuity of Gaussian volume, we get
$$\lim_{l\rightarrow \infty}\gamma_n(Q_l)=\gamma_n(K)=\frac{1}{2}.
$$

Now we prove the uniform lower bound of support function by contradiction. For $K_i\in\kno\rightarrow K$  with  $\gamma_n(K_i)= \frac{1}{2}$. Suppose that there exists $v_i\in\sn$ such that $h_{K_i}(v_i)\rightarrow 0$. Then for any $\epsilon>0,$  $v_i\in\sn$, $K_i\subset\left\{x|x\cdot v_i\geq-\epsilon\right\}$. Combining  the fact that $h_{K_i}(v_i)$ has upper bound, then  $K_i\subset B_R\cap\left\{x|x\cdot v_i\geq-\epsilon\right\}$. Let $H$ be halfspace, as computed in \cite{B},
$$
\gamma_n(\rn)=\frac{1}{(\sqrt{2\pi})^n}\int_{\rn}e^{\frac{-|x|^2}{2}}dx=\frac{1}{(\sqrt{2\pi})^n}\times\frac{1}{2}\times\sqrt{\frac{2}{\pi}}\times(\sqrt{2\pi})^{n+1}=1.
$$
Then
$$
\frac{1}{2}=\frac{1}{(\sqrt{2\pi})^n}\int_He^{\frac{-|x|^2}{2}}dx=\frac{1}{(\sqrt{2\pi})^n}\int_{B_{\frac{R}{2}}}e^{\frac{-|x|^2}{2}}dx+\frac{1}{(\sqrt{2\pi})^n}\int_{H\setminus B_{\frac{R}{2}}}e^{\frac{-|x|^2}{2}}dx,
$$
which shows that for some $t_0>0$,  $\frac{1}{(\sqrt{2\pi})^n}\int_{B_{\frac{R}{2}}}e^{\frac{-|x|^2}{2}}dx+t_0<\frac{1}{2}$, and naturally for $\epsilon$ small enough, there is $\gamma_n(K_i)\leq\gamma_n(B_R\cap\left\{x|x\cdot v_i\geq-\epsilon\right\})<\frac{1}{2}$, this contradicts to the condition $\gamma_n(K_i)=\frac{1}{2}$.
 Therefore, we conclude that $K$ is non-degenerate. Namely, $K$ is the desired convex body.
\end{proof}

Now we return to deal with the origin problem.
\subsection{Existence of  solution to the $L_p$-Gaussian Minkowski problem}

Combining with Lemma \ref{If}, Lemma \ref{var} and Theorem \ref{Eoo}, we prove sufficiently the main existence theorem for the $L_p$-Gaussian Minkowski problem stated in the introduction.

\begin{thm}
{\it For  $p>0,$ let $\mu$ be a nonzero finite  Borel measure on $\sn$ and be not concentrated in any closed hemisphere. There exist $K\in\kno$  and constant $\lambda>0$ such that
\begin{equation*}
\mu=\frac{\lambda}{p} S_{p,\gamma_n,K}.
\end{equation*}
 }
\end{thm}

The existence of a solution has been done as above, now we are devoted to solving  the  uniqueness.

\section{Uniqueness of solution}\label{SecU1}
In general, Brunn-Minkowski inequality implies the uniqueness of solution of Minkowski problem. There are many fruitful results on the inequality of Gaussian volume $\gamma_n$, such as  \cite{BC,EM,GZ}. In particular, Ehrhard inequality  in \cite{E} is one of the most significant Brunn-Minkowski type inequalities for Gaussian measure $\gamma_n$,  stated as follows. But in Gaussian Minkowski problem, because there is no homogeneity, the uniqueness can not obtained from the Brunn-Minkowski inequality. Fortunately, by using Ehrhard inequality, we can get the uniqueness result when the Gaussian volume is more than or equal to half.
\begin{thm}[\bf Ehrhard inequality]Let $K,L$ be convex bodies in $\rn$, $0\leq\lambda\leq1$, then
$$
\Phi^{-1}(\gamma_n(\lambda L+(1-\lambda)K))\geq\lambda\Phi^{-1}(\gamma_n(L))+(1-\lambda)\Phi^{-1}(\gamma_n(K)),
$$
with equality holds if and only if $K=L.$ Where $\Phi(x)=\frac{1}{\sqrt{2\pi}}\int_{-\infty}^xe^{\frac{-t^2}{2}}dt$.
\end{thm}
 In view of Ehrhard inequality, we obtain  the log-concave property of  $\gamma_n$ as below.
\begin{lem}
 Let $K,L$ be  convex bodies in $\rn$,  $0\leq\lambda\leq1$, then
 \begin{equation}
 \gamma_n(\lambda L+(1-\lambda) K)\geq {\gamma_n(L)}^\lambda{\gamma_n(K)}^{1-\lambda}
 \end{equation}
  with equality holds if and only if $K=L.$
\end{lem}
 Indeed,  this follows directly from  the above Lemma and  the fact that $\lambda L+_p(1-\lambda)\cdot K\supset\lambda L+_{p'}(1-\lambda)\cdot K$
for $p>p'$.
\begin{lem}
 Let $K,L$ be  convex bodies in $\rn$, for $p\geq 1$ and $0\leq\lambda\leq1$, then
 \begin{equation}
 \gamma_n(\lambda L+_p(1-\lambda)\cdot K)\geq {\gamma_n(L)}^\lambda{\gamma_n(K)}^{1-\lambda}
 \end{equation}
  with equality holds if and only if $K=L.$ And its differential equivalent form is
  \begin{equation}
  \frac{1}{p}\int_{\sn}{h_L}^pdS_{p,\gamma_n,K}\geq \frac{1}{p}\int_{\sn}{h_K}^pdS_{p,\gamma_n,K}+\gamma_n(K)\log\frac{\gamma_n(L)}{\gamma_n(K)}.
  \end{equation}
\end{lem}

Naturally, it tells that
\begin{lem}\label{15}
Let $K,L$ be convex bodies in $\rn$. For $p\geq 1$, if $\gamma_n(K)=\gamma_n(L),$ then
\begin{equation*}
\int_{\sn}{h_L}^pdS_{p,\gamma_n,K}\geq\int_{\sn}{h_K}^pdS_{p,\gamma_n,K}
\end{equation*}
with equality holds if and only if $K=L.$
\end{lem}

Finally, we attempt to deal with the uniqueness. The main idea of proof is inspired by Lemma 5.1 in \cite{HXZ}.
\begin{lem}\label{volume}
 For $p\geq1$, if $K,L \in \kno$  with the same $L_p$-Gaussian surface area, i.e.,
\begin{equation*}
S_{p,\gamma_n,K}=S_{p,\gamma_n,L}.
\end{equation*}
If $\gamma_n(K), \gamma_n(L)\geq\frac{1}{2},$  then $\gamma_n(K)=\gamma_n(L).$
\end{lem}
\begin{proof}
By the Ehrhard inequality,
\begin{equation*}
\Phi^{-1}\left(\gamma_n((1-t)K+tL)\right)\geq (1-t)\Phi^{-1}(\gamma_n(K))+t\Phi^{-1}(\gamma_n(L))
\end{equation*}
with equality holds if and only if $K=L.$ Where
\begin{equation*}
\Phi(x)=\frac{1}{\sqrt{2\pi}}\int_{-\infty}^x e^{\frac{-t^2}{2}}dt.
\end{equation*}

For convenience, we write $\Psi=\Phi^{-1}$, then
\begin{equation*}
\Psi(\gamma_n((1-t)K+tL))\geq  (1-t)\Psi(\gamma_n(K))+t\Psi(\gamma_n(L)).
\end{equation*}

For $p\geq1,$ $(1-t)K+tL\subset (1-t)\cdot K+_pt\cdot L$, together with the fact that $\Psi$ is $C^\infty$ and strictly monotonically increasing, it follows  that
\begin{equation*}
\Psi\left(\gamma_n\left((1-t)\cdot K+_pt\cdot L\right)\right)\geq \Psi(\gamma_n((1-t)K+tL))\geq (1-t)\Psi(\gamma_n(K))+t\Psi(\gamma_n(L)).
\end{equation*}
Then  by applying the Theorem  \ref{Variation} we have,
\begin{equation}\label{inter}
\frac{\Psi'\left(\gamma_n(K)\right)}{p}\int_{\sn}h_L^p-h_K^pdS_{p,\gamma_n,K}\geq \Psi(\gamma_n(L))-\Psi(\gamma_n(K)).
\end{equation}
Interchanging the position of $K$ and $L$, we have
\begin{equation*}
\frac{\Psi'\left(\gamma_n(L)\right)}{p}\int_{\sn}h_K^p-h_L^pdS_{p,\gamma_n,L}\geq \Psi(\gamma_n(K))-\Psi(\gamma_n(L)),
\end{equation*}
or equivalently,
\begin{equation}\label{p}
\frac{\Psi'\left(\gamma_n(L)\right)}{p}\int_{\sn}h_L^p-h_K^pdS_{p,\gamma_n,L}\leq \Psi(\gamma_n(L))-\Psi(\gamma_n(K)).
\end{equation}

By the condition $S_{p,\gamma_n,K}=S_{p,\gamma_n,L}$, $\Psi'(x)=\sqrt{2\pi}e^{\frac{\Psi(x)^2}{2}}>0$, \eqref{inter} and \eqref{p}, we have
\begin{equation*}
\frac{\Psi(\gamma_n(L))-\Psi(\gamma_n(K))}{\Psi'(\gamma_n(L))}\geq \frac{\Psi(\gamma_n(L))-\Psi(\gamma_n(K))}{\Psi'(\gamma_n(K))},
\end{equation*}
i.e.,
\begin{equation}\label{ineq}
\left(\Psi'(\gamma_n(K))-\Psi'(\gamma_n(L))\right)\left(\Psi(\gamma_n(K))-\Psi(\gamma_n(L))\right)\leq 0.
\end{equation}

On the one hand, $\Psi''(x)=\sqrt{2\pi}e^{\frac{\Psi(x)^2}{2}}\Psi(x)\Psi'(x)>0$ on $[\frac{1}{2}, \infty]$, $\Psi'$ is strictly increasing. On the other hand,  $\Psi$ is also strictly increasing, then we conclude
 \begin{equation}\label{Ineq}
 \left(\Psi'(\gamma_n(K))-\Psi'(\gamma_n(L))\right)\left(\Psi(\gamma_n(K))-\Psi(\gamma_n(L))\right)\geq 0
 \end{equation}
 with equality holds if and only if $\gamma_n(K)=\gamma_n(L)$. Combining \eqref{ineq} and \eqref{Ineq}, then $\gamma_n(K)=\gamma_n(L).$
\end{proof}
\begin{thm}\label{u}
 For $p\geq1$, if $K,L \in \kno$  with the same $L_p$-Gaussian surface area, i.e.,
\begin{equation*}
S_{p,\gamma_n,K}=S_{p,\gamma_n,L}.
\end{equation*}
If $\gamma_n(K), \gamma_n(L)\geq\frac{1}{2},$  then $K=L.$
\end{thm}
\begin{proof}
From Lemma \ref{volume}, we get $\gamma_n(K)=\gamma_n(L).$   Combining  the result of Lemma \ref{15} with the condition $S_{p,\gamma_n,K}=S_{p,\gamma_n,L}$, it sufficient to have
\begin{equation*}
\int_{\sn}{h_L}^pdS_{p,\gamma_n,L}\geq\int_{\sn}{h_K}^pdS_{p,\gamma_n,K},
\end{equation*}
 by changing the position of $K,L$,  the equality holds in Lemma \ref{15}, i.e., $K=L.$
\end{proof}

The following isoperimetric type inequality can be derived from Ehrhard inequality, and will be used in the  Section~\ref{SecE}.
\begin{thm}\label{13}
For $p\geq1,$ let $K$ be an $o$-symmetry convex body in $\rn,$ and symmetric strip $P=\left\{|x_1|\leq a\right\}$ with $\gamma_n(K)=\gamma_n(P)=\frac{1}{2}$. Then,
$$
|S_{p,\gamma_n,K}|\geq \sqrt{\frac{2}{\pi}}r^{-p}ae^{\frac{-a^2}{2}},
$$
where $r$ is chosen such that $\gamma_n(rB)=\frac{1}{2}$.
\end{thm}
\begin{proof}
For $p\geq1,$ applying Ehrhard inequality, $\Phi^{-1}$ is monotonically increasing, and
$(1-t)K+tL\subset (1-t)\cdot K+_pt\cdot L$,
then,
\begin{equation*}
\begin{split}
\Phi^{-1}(\gamma_n((1-t)\cdot K+_pt\cdot L))&\geq\Phi^{-1}\left(\gamma_n((1-t)K+tL)\right)\\
&\geq(1-t) \Phi^{-1}(\gamma_n(K))+t\Phi^{-1}(\gamma_n(L)).
\end{split}
\end{equation*}

Set $L=rB$ such that $\gamma_n(rB)=\gamma_n(K)=\frac{1}{2}$, it follows that
\begin{equation}\label{low}
|S_{p,\gamma_n,K}|\geq r^{-p}\int_{\sn}h_KdS_{\gamma_n,K}.
\end{equation}

By the known result in \cite{LR},
$$
\gamma_n(K)=\gamma_n(P)\Rightarrow  \left. \frac{d}{d t} \right|_{t=1}\gamma_n(tK)\geq\left. \frac{d}{d t} \right|_{t=1}\gamma_n(tP),
$$
i.e.,
\begin{equation}\label{lower}
\int_{\sn}h_KdS_{\gamma_n,K}\geq \sqrt{\frac{2}{\pi}}ae^{\frac{-a^2}{2}}.
\end{equation}

Combining the \eqref{low} and  \eqref{lower}, we have
$$
|S_{p,\gamma_n,K}|\geq \sqrt{\frac{2}{\pi}}r^{-p}ae^{\frac{-a^2}{2}}.
$$
\end{proof}

\section{Existence of weak solution}\label{SecE}
In this section, we do some prior estimates and apply degree theory to  obtain the existence of  $o$-symmetry  smooth solution, then with the help of approximation argument to get the existence of weak solution of $L_p$-Gaussian Minkowski problem for $p\geq1$.
\subsection{prior estimates}
\begin{lem}[\bf $C^0$ estimate]\label{19}
For $p\in\rbo$, $K\in\kno$ and $\gamma_n(K)\geq \frac{1}{2}$, assume $h_K\in C^2(\sn)$ is the solution of the following equation
\begin{equation*}
\frac{1}{(\sqrt{2\pi})^n}e^{-\frac{|\nabla h|^2+h^2}{2}}h^{1-p}\det(\nabla^2h+hI)=f.
 \end{equation*}
 If $\frac{1}{C_1}<f<C_1$ for some positive constant $C_1$, then there exists a constant  $C_2>0$ such that $\frac{1}{C_2}<h_K<C_2$.
\end{lem}
\begin{proof}
Now we prove the upper bounded. There exists $v_0\in\sn$ such that $h_K(v_0)=\max_{v\in\sn}h_K(v)$,  $\nabla h_K|_{v_0}=0$ and $ \nabla^2h_K|_{v_0}\leq 0$, then
\begin{align*}
f(v_0)&=\frac{1}{(\sqrt{2\pi})^n}[e^{-\frac{|\nabla h|^2+h^2}{2}}h^{1-p}\det(\nabla^2h+hI)](v_0)\\
&\leq \frac{1}{(\sqrt{2\pi})^n}e^\frac{-h_K(v_0)^2}{2}h_K(v_0)^{n-p},
\end{align*}
which implies $h_K(v_0)\leq C_2.$

 The bound from blew of support function is guaranteed  by the condition that $\gamma_n(K)\geq\frac{1}{2}$. The  proof  is same as established as in Theorem \ref{Eoo}.
\end{proof}

\begin{lem}
For $p\geq1,$ $\alpha\in(0,1)$, assume function $f\in C^{2,\alpha}(\sn)$ and there exists  constant $C>0$ such that $|f|_{C_{2,\alpha}}<C$ and $\frac{1}{C}<f<C$. For  $K\in\kno$ with $\gamma_n(K)>\frac{1}{2}$, if $h_K\in C^{4,\alpha}(\sn)$ and satisfies
\begin{equation}\label{MAE}
\frac{1}{(\sqrt{2\pi})^n}h^{1-p}e^{\frac{-(|\nabla h|^2+h^2)}{2}}\det(\nabla^2h+hI)=f,
\end{equation}
then  there exists a positive constant $C'$ which only depends on $C$ such that the following  priori estimates hold:

(1)  $C^1$ estimate: $|\nabla h_K|<C'.$

(2)   $C^2$ estimate: $\frac{1}{C'}I<(\nabla^2h_K+h_KI)<C'I$ and  higher estimate $|h_K|_{C^{4,\alpha}}<C'$.
\end{lem}
\begin{proof}
Assume  $M(v_0)=\max_{v\in\sn}(|\nabla h_K(v)|^2+h_K(v)^2)$.  It admits
$$
2h_Kh_{Ki}+2h_{Kl}h_{Kli}=0,
$$
at $v_0$, i.e.,
$$h_{Kl}(h_K\delta_{il}+h_{Kli})=0.$$
In view of \eqref{MAE}, it tells that  the element of the matrix $h_K\delta_{il}+h_{Kli}>0,$ we have $\nabla h_K(v_0)=0$,  and $M(v_0)=h_K(v_0)^2<C'$, the second inequality is immediately  from Lemma \ref{19} and $C'$ is only depends on $C$.

To get $C^2$ estimate,  the key is to show that the eigenvalues of matrix $\nabla^2h+hI$  are bounded from above and below. In other words,  the Monge-Ampere equation \eqref{MAE} is uniformly elliptic. Thus, an immediate consequence of the standard Evans-Krylov-Safonov theory \cite{GT} is $|h_K|_{C^{4,\alpha}}<C'$.

On the one hand, in light of the fact that $f, h_K$ have positive upper and lower bounds,  then combined with $|\nabla h_K|<C'$, we conclude
$$
\det(\nabla^2h_K+h_KI)=(\sqrt{2\pi})^n{h_K}^{p-1}fe^{\frac{|\nabla h_K|^2+{h_K}^2}{2}}
$$
also has positive upper and lower bounds.

On the other hand, we end the proof by claiming that the trace of $\nabla^2h_K+h_KI$ has an upper bound.

 For convenience, denote
 \begin{equation*}
 H=trace(\nabla^2h_K+h_KI)=\triangle h_K+(n-1)h_K.
 \end{equation*}

 Assume $H(v_1)=\max_{v\in\sn}H(v)$, then $\nabla H(v_1)=0$ and $\nabla^2H(v_1)\leq 0$. We can make the Hessian of $h_K$, $(h_K)_{ij}$, is diagonal by choosing the suitable local orthogonal frame $e_i(i=1,2,\ldots,n)$.  Denote $\omega_{ij}=(h_K)_{ij}+h_K\delta_{ij}$, setting $\omega^{ij}=\omega_{ij}^{-1}$,  the inverse matrix of $\omega_{ij}$.

 Then,  based on the above conclusions, at $v_1,$ we have
 \begin{equation}\label{0}
 0\geq \omega^{ij}H_{ij}=\omega^{ii}H_{ii}=\omega^{ii}\triangle\omega_{ii}-(n-1)^2+H\sum_{i=1}^{n-1}\omega^{ii}\geq \omega^{ii}\triangle\omega_{ii}-(n-1)^2,
 \end{equation}
 the second equality of \eqref{0} is from the commutator identity \cite{GMZ}:
 \begin{equation*}
 H_{ii}=\triangle\omega_{ii}-(n-1)\omega_{ii}+H.
 \end{equation*}

 Next we are going to  estimate $\omega^{ii}\triangle\omega_{ii}.$

  From the equation
 \begin{equation*}
 \frac{1}{(\sqrt{2\pi})^n}h_K^{1-p}e^{\frac{-(|\nabla h_K|^2+h_K^2)}{2}}\det(\nabla^2h_K+h_KI)=f,
 \end{equation*}
i.e.,
 \begin{equation*}
 \log \det(\nabla^2h_K+h_KI)=\log f+\frac{|\nabla h_K|^2+h_K^2}{2}+(p-1)\log h_K+\frac{n}{2}\log (2\pi).
 \end{equation*}

 Carry out differential operation twice on the above formula, we obtain
\begin{equation*}
\omega^{ij}\omega_{ij\alpha}=(\log f)_\alpha+h_K(h_K)_\alpha+(h_K)_l(h_K)_{l\alpha}+(p-1)\frac{(h_K)_\alpha}{h_K},
 \end{equation*}
 and

 \begin{equation*}
 \begin{split}
 &\omega^{ij}\omega_{ij\alpha\alpha}+(\omega^{ij})_\alpha(\omega_{ij\alpha})\\
 &=(\log f)_{\alpha\alpha}+h_K(h_K)_{\alpha\alpha}+((h_K)_\alpha)^2+(h_K)_{l\alpha}^2
 +(h_K)_l(h_K)_{l\alpha\alpha}\\
 &+(p-1)\frac{h_K(h_K)_{\alpha\alpha}-(h_K)_\alpha^2 }{h_K^2},
 \end{split}
 \end{equation*}
 it illustrates that
  \begin{equation}\label{21}
 \begin{split}
&\omega^{ij}\omega_{ij\alpha\alpha}= -(\omega^{ij})_\alpha(\omega_{ij\alpha})+(\log f)_{\alpha\alpha}+h_K(h_K)_{\alpha\alpha}+((h_K)_\alpha)^2+(h_K)_{l\alpha}^2\\
 &+(h_K)_l(h_K)_{l\alpha\alpha}
 +(p-1)\frac{h_K(h_K)_{\alpha\alpha}-(h_K)_\alpha^2 }{h_K^2}.
\end{split}
 \end{equation}
 We estimate the right side of the \eqref{21}.

Since
 \begin{equation*}
 (\omega^{ij}\omega_{ij})_\alpha=(1)_\alpha=0,
 \end{equation*}
 i.e.,
  \begin{equation*}
  (\omega^{ij})_\alpha\omega_{ij}+\omega^{ij}\omega_{ij\alpha}=0,
  \end{equation*}
then we have
\begin{equation*}
(\omega^{ij})_\alpha=-\omega^{ij}\omega^{ij}\omega_{ij\alpha},
\end{equation*}
and
\begin{equation*}
(\omega^{ij})_\alpha(\omega_{ij\alpha})=-(\omega^{ij})^2(\omega_{ij\alpha})^2\leq 0.
\end{equation*}
From the equation (4.11) in Cheng-Yau\cite{CY}, we have
\begin{equation*}
 \begin{split}
&\sum_i(h_K)_i(\Delta (h_K)_i)\\
&=\sum_i[(h_K)_i(\Delta h_K)_i+(h_K)_i(h_K)_i-(h_K)_i(h_K)_\alpha\delta_{i\alpha}]\\
&=\nabla h_K\cdot\nabla H-(n-1)|\nabla h_K|^2.
 \end{split}
\end{equation*}
Due to  $H=\Delta h_K+(n-1)h_K$, $\nabla H(v_1)=0$, $\sum_{\alpha=1}^{n-1}h_{\alpha\alpha}^2\geq \frac{(\sum_{\alpha=1}^{n-1} h_{\alpha\alpha})^2}{2}=\frac{(H-(n-1)h_K)^2}{2}$, thus, at $v_1$,
 \begin{equation}\label{1}
 \begin{split}
 &\sum_{\alpha=1}^{n-1}\omega^{ii}\omega_{ii\alpha\alpha}\\
 &\geq\Delta\log f+h_K(H-(n-1)h_K)+\frac{(H-(n-1)h_K)^2}{2}\\
 &-(n-2)|\nabla h_K|^2+(p-1)\frac{H}{h_K}-(p-1)(n-1)+(1-p)\frac{|\nabla h_K|^2}{h_K^2}\\
 &=\frac{H^2}{2}+\left[(2-n)h_K+\frac{p-1}{h_K}\right]H+\Delta\log f-\frac{(n-1)^2h_K^2}{2}\\
 &-(n-2)|\nabla h_K|^2-(p-1)(n-1)-(p-1)\frac{|\nabla h_K|^2}{h_K^2}.
 \end{split}
 \end{equation}
 Put this in \eqref{0}, then
  \begin{equation}
 \begin{split}
 0&\geq\frac{H^2}{2}+\left[(2-n)h_K+\frac{p-1}{h_K}\right]H+\Delta\log f-\frac{(n-1)^2h_K^2}{2}\\
 &-(n-2)|\nabla h_K|^2-(p-1)(n-1)-(p-1)\frac{|\nabla h_K|^2}{h_K^2}-(n-1)^2.
 \end{split}
 \end{equation}
Since $f$, $|f|_{C^{2,\alpha}}$, $h_K$, $|\nabla h_K|$ are bounded, which implies  $H\leq C'.$

Next we focus on  obtaining the existence of smooth solution by  the degree theory for second-order nonlinear elliptic operators, the reader can conference to the Li \cite{L} for some details.
\end{proof}
\begin{thm}[\bf Existence and uniqueness of smooth solution]\label{14}
Let $\alpha\in (0,1)$,  for $p\geq 1$, $f\in C^{2,\alpha}(\sn)$ which is positive even function and satisfies $|f|_{L_1}<\sqrt{\frac{2}{\pi}}r^{-p}ae^{\frac{-a^2}{2}}$, the $r$ and $a$ are chosen such that $\gamma_n(rB)=\gamma_n(P)=\frac{1}{2}$, symmetry trip $P=\left\{x\in \rn:|x_1|\leq a\right\}$, then there exists a unique $C^{4,\alpha}$ convex body  $K\in\kne$ with $\gamma_n(K)>\frac{1}{2}$, its support function $h_K$ satisfies
\begin{equation}\label{12}
\frac{1}{(\sqrt{2\pi})^n}e^{-\frac{|\nabla h_K|^2+h_K^2}{2}}h_K^{1-p}\det(\nabla^2h_K+h_KI)=f.
\end{equation}
\end{thm}
\begin{proof}
 On the one hand,  the uniqueness result is guaranteed by  Theorem \ref{u}  .

On the other hand,  define  $F(;t):C^{4,\alpha}(\sn)\rightarrow C^{2,\alpha}(\sn)$ as follows:
\begin{equation*}
	F(h;t) = \det(\nabla^2 h+hI)- {(\sqrt{2\pi})^n} e^{\frac{|\nabla h|^2+h^2}{2}}h^{p-1}f_t,
\end{equation*}	
where $f_t = (1-t)c_0+tf$.

Define $O\subset C^{4,\alpha}(\sn)$ by
\begin{align*}
	O=\left\{h: \frac{1}{C'}<h<C', \frac{1}{C'}I<(\nabla^2h + hI)<C'I, \right.\left.F(h;t)=0, |h|_{C^{4,\alpha}}<C', \gamma_n(h)>\frac{1}{2}\right\}.
\end{align*}
For $h\in O$, the eigenvalues of its hessian are bounded from above and below,  the operator $F(\cdot;t)$ is uniformly elliptic on $O$ for any $t\in[0,1].$

When $f=c_0>0(t=0)$ is small enough, applying mean value theorem, it reveals  that there exists a unique constant solution $h_K=r_0$ such that $\gamma_n(K)>\frac{1}{2}$. Since spherical Laplacian has a discrete spectrum, then we can  select $c_0$ suitably such that  $|c_0|_{L_1}<\sqrt{\frac{2}{\pi}}r^{-p}ae^{\frac{-a^2}{2}}$, and the operator $L\phi=\Delta_{\sn}\phi+((n-p)-{r_0}^2)\phi$ is invertible.

Now we claim that $O$ is an open bounded set under the norm $|\cdot|_{C^{4,\alpha}}$, that is, we need to prove that if $h\in\partial O$, then $F(h;t)\neq0$.

If $F(h;t)=0$, in other words, $h$ is the solution of
$$
 \frac{1}{(\sqrt{2\pi})^n} e^{-\frac{|\nabla h|^2+h^2}{2}}h^{1-p}\det(\nabla^2 h+hI)=f_t.
$$
Since $h\in\partial O,$ then $\gamma_n([h])=\frac{1}{2}$, from Theorem \ref{13}, we have
$|S_{p,\gamma_n,[h]}|\geq \sqrt{\frac{2}{\pi}}r^{-p}ae^{\frac{-a^2}{2}}$, which is contracted  to the condition $|f_t|_{L_1}<\sqrt{\frac{2}{\pi}}r^{-p}ae^{\frac{-a^2}{2}}$.

By means of the Proposition 2.2 of Li\cite{L}, we  conclude that
$$
deg(F(\cdot;0),O,0)=deg(F(\cdot;1),O,0).
$$
It is clear that  if $deg(F(\cdot;1),O,0)\neq0,$ then there exists $h\in O$ such that $F(h;1)=0$. Subsequently, we need to claim   $deg(F(\cdot;0),O,0)\neq0.$ The Proposition 2.3 and Proposition 2.4  of Li \cite{L} told that if $L_{r_0}$ the linearized operator of $F$ at $r_0$ is invertible, then we have
$$
deg(F(\cdot;0),O,0)=deg(L_{r_0}(\cdot;0),O,0)\neq0.
$$

Our final goal is to verify  $L_{r_0}$ :$ C^{4,\alpha}(\sn)\rightarrow C^{2,\alpha}(\sn)$,
\begin{equation*}
\begin{split}
L_{r_0}&={r_0}^{n-2}\triangle_{\sn}\phi+((n-p){r_0}^{n-2}-{r_0}^{n})\phi\\
&={r_0}^{n-2}(\triangle_{\sn}\phi+((n-p)-{r_0}^{2})\phi)
\end{split}
\end{equation*}
is invertible.

The choice of $c_0$ ensures the reversibility of  $L_{r_0}$, then we completed  the claim.
\end{proof}

We via approximation argument to obtain the existence of weak solution of $L_p$-Gaussian Minkowski problem as follows.
\begin{lem}\label{18}
 If $\mu$ is not concentrated on any closed hemisphere, let $\mu_i=f_idv$ be a sequence of measure which  converges to $\mu$ weakly as $i\rightarrow\infty.$  Then there exist  $\epsilon_0>0, \delta_0>0, $ and $N_0>0$ such that for any $i>N_0, e\in\sn$,
 $$
 \int_{\sn\cap\left\{v|v\cdot e>\delta_0\right\}}f_idv>\epsilon_0.
 $$
\end{lem}
\begin{proof}
Argue via contradiction.  Take $\epsilon_0=\frac{1}{k},\delta_0=\frac{1}{k}$, $e_k\in\sn$, a subsequence which converges to $e,$
 we conclude that
$$
 \int_{\sn\cap\left\{v|v\cdot e_k>\frac{1}{k}\right\}}f_kdv\leq\frac{1}{k}.
 $$

Since  $f_idv\rightharpoonup \mu$ weakly, fix any $\eta>0,$ for $k$ large enough,
 $$
 \int_{\sn\cap\left\{v|v\cdot e>\eta\right\}}d\mu= 0.
 $$
 Let $\eta\rightarrow 0,$ we have $\int_{\sn\cap\left\{v|v\cdot e>0\right\}}d\mu= 0$, which is contradicted to the condition that $\mu$ is not concentrated on any closed hemisphere.
\end{proof}

At the end, we are ready to  state the Theorem for the existence of unique symmetry weak solution as exhibited in the following.
\begin{thm}\label{exi}
For $p\geq1$, let $\mu$ be a finite even Borel measure on $\sn$ and be not concentrated in any closed hemisphere with $|\mu|<\sqrt{\frac{2}{\pi}}r^{-p}ae^{\frac{-a^2}{2}}$, the $r$ and $a$ are chosen such that $\gamma_n(rB)=\gamma_n(P)=\frac{1}{2}$, symmetry trip $P=\left\{x\in \rn:|x_1|\leq a\right\}$. Then there exists a unique $K\in\kne$ with $\gamma_n(K)>\frac{1}{2}$ such that
$$
S_{p,\gamma_n,K}=\mu.
$$
\end{thm}
\begin{proof}

Let $\mu_i=f_idv$ be a sequence of even measure which  converges to $\mu$ weakly as $i\rightarrow\infty,$  $f_i\in C^{2,\alpha}, f_i>0$ with $0<|f_i|_{L_1}<\sqrt{\frac{2}{\pi}}r^{-p}ae^{\frac{-a^2}{2}}$. Then combine Theorem \ref{14}, there are $C^{4,\alpha}$ convex bodies $ K_i\in\kne$ with $\gamma_n(K_i)>\frac{1}{2}$, its support function satisfies
$$
 \frac{1}{(\sqrt{2\pi})^n} e^{-\frac{|\nabla h_{K_i}|^2+{h_{K_i}}^2}{2}}{h_{K_i}}^{1-p}\det(\nabla^2 h_{K_i}+h_{K_i}I)=f_i,
$$
i.e., $S_{p,\gamma_n,K_i}=f_idv.$

We will set about claiming that $K_i$ is uniformly bounded. Suppose  that $\max_{v\in\sn}h_{K_i}(v)$ is attained at $v_i\in\sn$, and the corresponding point on $K$ is $x_i.$ In fact, $|x_i|=h_{K_i}(v_i)$. Then  by Lemma \ref{18}, there exist $\epsilon_0>0, \delta_0>0,$ and $N_0>0$ such that for any $i>N_0, $
 \begin{equation}\label{20}
 \int_{\sn\cap\left\{v|v\cdot v_i>\delta_0\right\}}f_idv>\epsilon_0.
\end{equation}

Denote $D_i=\sn\cap\left\{v|v\cdot v_i>\delta_0\right\}$, the points $x$ on $K_i$ corresponding to the vector on $D_i$ satisfy $|x|>h_{K_i}(v_i)\delta_0$. Hence,
\begin{equation*}
\begin{split}
S_{p,\gamma_n,D_i}&=\frac{1}{(\sqrt{2\pi})^n}\int_{\nu_{K_i}^{-1}(D_i)} e^{-\frac{|x|^2}{2}}(x\cdot\nu_{K_i}(x))^{1-p}d\hm(x)\\
&\leq\frac{1}{(\sqrt{2\pi})^n}e^{-(h_{K_i}(v_i)\delta_0)^2}\left(h_{K_i}(v_i)\delta_0\right)^{1-p}{h_{K_i}(v_i)}^{n-1}\rightarrow 0,
\end{split}
\end{equation*}

as $h_{K_i}(v_i)\rightarrow\infty.$ Which  contradicts to \eqref{20}.

In view of the fact that $K_i$ is uniformly bounded, $\gamma_n(K_i)>\frac{1}{2}$,  it follows that $h_{K_i}>\frac{1}{C}$.
In summary, we have $\frac{1}{C}<h_{K_i}<C,$   $K_i\rightarrow K\in\kne$ as  $i\rightarrow\infty.$ $K$ is the desired convex body.
\end{proof}

\section*{Acknowledgement}
 I would like to thank  my supervisor, professor Yong Huang, for his patient guidance and encouragement.   I am also deeply indebted to professor Shibing Chen for providing the valuable advice of Lemma~\ref{19}.

\end{document}